\documentclass[11pt]{article}\hfuzz=10pt  \topmargin=-0.5cm\hfuzz=10pt  \oddsidemargin=0.1cm\textheight=220mm \textwidth=16.5cm
\usepackage{amssymb}
\usepackage{epsfig}
\newcommand{\comm}[1]{ }

\def\xxxonly{\comm}
\def\xxxonly{{ }}
\def\noxxx{\comm}

\usepackage[pagewise]{lineno}
 \usepackage{color}

\newcommand{\bl}[1]{\color{blue}#1\color{black}}

\def\bl{\comm}

\newtheorem{theorem}{Theorem}
\newtheorem{lemma}{Lemma}
\newtheorem{proposition}{Proposition}
\newtheorem{remark}{Remark}
  
\def\defi{\stackrel{{\scriptscriptstyle \Delta}}{=}}

\def\a{\alpha}
\def\d{\delta}
\def\o{\omega}

\def\w{\widehat}

\def\Re{{\rm Re\,}}

\def\R{{\bf R}}

\def\H{{\cal H}}

\def\L{L}

\def\b{\beta}

\def\g{\gamma}
\def\C{{\bf C}}
\def\W{{\cal W}^*}

\def\t{\theta}
\def\oo{\bar}

\def\L{{\cal L}}

\newcommand{\be}{\begin{equation}}
\newcommand{\ee}{\end{equation}}
\newcommand{\bd}{\begin{displaymath}}
\newcommand{\ed}{\end{displaymath}}
\newcommand{\ba}{\begin{array}{ll}}
\newcommand{\ea}{\end{array}}
\newcommand{\baa}{\begin{eqnarray}}
\newcommand{\eaa}{\end{eqnarray}}
\newcommand{\baaa}{\begin{eqnarray*}}
\newcommand{\eaaa}{\end{eqnarray*}}

\def\L{{\cal L}}

\def\t{\tau}

\def\W{{\cal W}}

\def\DD{\oo{\cal D}}
\def\o{\omega}
\def\H{{\cal H}}

\def\R{{\bf R}}

\def\t{\theta}

\def\H{H_{BC}}
\def\H{{H^1}}
\def\HH{{H^2}}
\def\HBC{{H_{BC}}}
\def\Im{{\rm Im}\,}

    \xxxonly{\date{Submitted: 8 July  2019. Revised: February 9 2021} 
\title{Regularity of complexified hyperbolic  equations with integral conditions}

\author{
Nikolai Dokuchaev}
\begin{document}
\maketitle
\let\thefootnote\relax
\noxxx{\footnote{
 }
 \let\thefootnote\relax\footnote{
The author is with 
Zhejiang University/University of Illinois at Urbana-Champaign Institute,  Zhejiang University, Haining, Zhejiang, China. Email: Dokuchaev@intl.zju.edu.cn}}
\vspace{-1cm}
\begin{abstract}
This  paper  considers  hyperbolic wave equations with non-local in time
conditions involving integrals with respect to time.
It is shown that regularity of the  solution can be achieved for complexified  problem with
integral conditions involving harmonic complex  exponential weights.
  The paper establishes existence, uniqueness,
 and  a regularity of the solutions.
\end{abstract}
\section{Introduction}
The most common type of boundary conditions  for
evolution partial differential equations are Cauchy initial conditions.
It is know that these conditions can be replaced, in some cases, by non-local conditions, for example, including integrals over time intervals.
There is  a significant number of works devoted to these boundary values problems.

For parabolic equations with non-local in time boundary conditions, some results and references    can be found, e.g.,  in \cite{D19,Pri}.
For Schr\"odinger equations,  some results and references    can be found in, e.g., \cite{AS,BP,D19xxx,SM}.   
In \cite{AS,BP,SM},  the non-local in time condition was dominated by the initial value, and the approach was  based on the contraction mapping theorem. In  \cite{AS,SM}, the nonlocal conditions connected solutions in a finite set of times. In \cite{BP}, the conditions were quite general and  allowed to include integrals with respect to  time.
In \cite{D19xxx}, integral conditions without dominating initial value have been considered.

For hyperbolic equations with mixed highest order derivative,
regularity results were obtained  in \cite{Ass}.
In \cite{AB}, systems of hyperbolic equations with integral conditions have been considered.

For hyperbolic wave equations,  related results and the references   can be found, e.g., in \cite{Bou,Bel,DeM,Dm,KKK,KP,MB,PS}.
In \cite{Bou,Bel,DeM,Dm, MB}, regularity results  were obtained for hyperbolic wave equations with  a variety of integral conditions with respect to state variables.

 In \cite{KKK,KP,PS}, hyperbolic wave equations under  integral conditions with respect to
 time have been considered. In  \cite{KKK,KP}, the eigenfunction expansion method has been used, and the regularity result has been affected by the so-called "small denominators" ("small divisors") problem  that often causes
instability of solutions for hyperbolic wave equations with non-local in time conditions. The solvability was  obtained in \cite{KKK,KP} for the case where the spectrum for the inputs and solutions does not contain resonance points.
 In  \cite{PS},   a regularity condition  without these restrictions on the spectrum have been obtained for the hyperbolic wave equation with a Laplacian.
 This condition imposed certain restrictions on the kernel in the integral condition which has vanish with a certain rate at the end of  the time interval.

 The paper readdresses the problem of  regularity for solutions of  boundary value problems for hyperbolic wave equations with
non-local in time conditions. The paper suggests a complexification  of boundary value problems that ensures regularity of the solutions. This  complexification requires
to consider integral conditions
 $\int_0^Te^{i\o t}u(t)dt=g$, for a real nonzero $\o$ that can be arbitrarily small, where $u(t)$ is the solution of the hyperbolic wave equation with real coefficients that take values in an appropriate  Hilbert space.
This  allows to bypass, for this particular setting,  the "small denominators" (or "small divisors") problem.
We  establish existence, uniqueness,
 and  a regularity of the solutions for the complexified equation. The proofs are  based the spectral expansion, similarly to the setting from \cite{D19,GG,LV,D19xxx}.  The eigenfunction expansion for the solution is presented explicitly.
This allows  to derive a numerical solution.

 The rest of the paper  is organized as follows. In Section \ref{SecM}, we introduce
a boundary value problem with averaging over time,
and we  present the main result  (Theorem \ref{ThM}).
In Section \ref{SecP}, we present the proofs. Section \ref{SecN} gives  a 
numerical example of the impact of the presence of small $\o\neq 0$ on the appearance of small denominators..

 Section \ref{SecD} presents conclusions and discusses future research.
\subsection*{Some definitions}\label{SecS}
For a Banach space $X$,
we denote  the norm by $\|\cdot\|_{ X}$. For a Hilbert space $X$,
we denote  the inner product by $(\cdot,\cdot)_{ X}$. We denote
 the Lebesgue measure and
 the $\sigma $-algebra of Lebesgue sets in $\R^n$
by $\oo\ell _{n}$ and $ {\oo{\cal B}}_{n}$, respectively.

\def\DD{{\rm D}}
 Let $\DD\subset \R^n$ be a domain, and let $H=L_2(\DD,\oo{\cal B}_n, \oo\ell_{1}; \C)$ be the space of complex-valued functions.  Let $A:H\to H$ be a self-adjoint operator defined
 on an everywhere closed  subset $D(A)$ of $H$ such that if $u$ is a real valued function then $Au$ is also a real valued function.

 Let $\HBC$ be a set functions from $H$ such that $\HBC$ is everywhere dense in $H$, that the set $D(A)\cap H$ is everywhere dense in  $H$, and that $(u, Av)_H$ is finite and defined for $u,v\in \HBC$ as a continuous extension from $D(A)\times D(A)$.

Consider an eigenvalue problem \baa
 Av=-\lambda v,\quad v\in\HBC.  \label{EP}\eaa
We assume that this equation is satisfied for $v\in \HBC$  if
  \baaa
 (w,Av)_H=-\lambda (w,v,)_H \quad \forall w\in\HBC  \label{EPw}.\eaaa

Assume that  there exists
a  basis  $\{v_k\}_{k=1}^\infty\subset \HBC$  in $H$  such that
 $$(v_k,v_m)_{H}=0,\quad  k\neq m,\quad \|v_k\|_{H}=1,$$
and that $v_k$ are eigenfunctions for (\ref{EP}), i.e.,
\baa
 Av_k=-\lambda_k v_k, \label{EPk}\eaa
for some $\lambda_k\in (0,+\infty)$ such that $\lambda_k\to +\infty$ as $k\to +\infty$.

These assumptions imply that the operator $A$ is self-adjoint with respect to the boundary conditions  defined by the choice of $\HBC$.

For $q=-1,0,1,2$, let $H^q$  be  the Hilbert spaces obtained
as the closure of the set  $\HBC\cap D(A)$ in the  norms  \baaa
&&\|u\|_{H^q}\defi \left(\sum_{k=1}^\infty\lambda_k^q(u,v_k)_H^2\right)^{1/2}
 \eaaa
respectively. According to this definition, $H^0=H$.

Clearly, the bilinear  forms $(u,v)_H$,
$(v, Au)_H$ and  $(Av, Au)_H$ are well defined on   $H^{-1}\times \H$,  $\H\times \H$ and  $\HH\times \HH$, respectively, since they  can be extended continuously from $D(A)\times D(A)$ .

\par
For  $q=-1,0,1,2$, and $r\in[1,+\infty]$,   introduce the spaces
\baaa
&&{{\cal C}^0\defi C\left([0,T]; H\right),\qquad} {\cal C}^q\defi C\left([0,T]; H^q\right),
\qquad \L_r^{q}\defi L_{r}\bigl([ 0,T ],\oo{\cal B}_1, \oo\ell_{1};  H^q\bigr),
\eaaa
and the spaces
\baaa
 \W^1_r\defi\Bigl\{u\in {\cal C}^0\cap \L^{1}_r:\ \frac{du}{dt}\in {\cal C}^{0}\Bigl\}
\eaaa
considered as Banach spaces with the norms, respectively,
\baaa
\|u\|_{\W^1_r}\defi \|u\|_{ {\cal C}^0}+ \|u\|_{\L^1_r}+\left\|\frac{du}{d t}\right\|_{{\cal C}^0}.
\eaaa
\section{Problem setting and the main result}\label{SecM}

Let   $T>0$, $\o\in \R$, $\o\neq 0$, \index{$\varrho\in\C$,?}  $a\in\H$, and $g\in\HH$, be given. We consider the boundary value problem
 \baa
&&\frac{d^2 u}{d t^2}(t)=A u(t), \quad t\in (0,T),\label{eq}\\
&&u(0)=a,\label{a0}\label{pp}\\
&&\int_0^T e^{i\o t}u(t)dt=g.\label{ppp}
\eaa
For $u\in \W^{1}_1$, we accept that equation (\ref{eq}) is
satisfied  as an equality 
\baa
\frac{d u}{d t}(t)-\frac{d u}{d t}(s)=\int_s^t A u(r)dr\label{eq0}
\label{intA}\eaa
that holds in $H^{-1}$ for all $s,t$ such that $0\le s\le t\le T$, and that conditions (\ref{pp}) and (\ref{pp}) are
satisfied as equalities  in $H$.

\begin{theorem}
\label{ThM} Assume  that $e^{2i\o T}\neq 1$. In this case, for any  $a\in\H$ and $g\in \HH$, there exists a unique solution $u\in\W^{1}_\infty$.
 Moreover, there exists $c>0$  such
that
\baa \|u\|_{\W^{1}_\infty}\le c(\|a\|_\H+\|g\|_{\HH})\quad \forall a\in\H,\ g\in\HH.
\label{estp}
\eaa
Here $c>0$ depends only on $\HBC,A,T$, and $\o$. \end{theorem}

\par
By Theorem \ref{ThM}, problem
(\ref{eq})-(\ref{ppp}) is well-posed in the sense of Hadamard for $a\in\H$ and $g\in\HH$.
The proof of this theorem is given below; it is based on  explicit representation  of the solution $u$ for given $g$ and $a$ via  eigenfunction expansion.
It can be noted that, since $\o\ne 0$ and $\o\in\R$, the solution $u$  of problem (\ref{eq})-(\ref{ppp}) is not real valued even if both  $a$ and $g$ are real valued.

\begin{remark}{\rm   It is known that
real valued problem  (\ref{eq})-(\ref{ppp}) with $\o=0$ is solvable for some real valued $\a$ and $g$ but it
does not feature stable solutions due to "small denominator" problem arising for certain frequencies;  see, e.g. examples in \cite{PS}, p. 42.
For small  $\o\to 0$, the part $v$ of the solution $(v,w)$ from
Remark  \ref{rem1}    can be used as a stable approximation of the real valued solution of  problem (\ref{eq})-(\ref{ppp}) with $\o=0$. In this case,
$\sin(\o t)w(t)$ can be considered as some small stabilizing term.
}\end{remark}

\begin{remark}\label{rem1}{\rm
One can reformulate the  setting with complex  valued solutions as a  setting with
two real valued solution. Assume that $u(t)=v(t)+iw(t)$, where $v(t)=\Re u(t)$ and $w(t)=\Im u(t)$. Then problem (\ref{eq})-(\ref{ppp})  can be rewritten as
\baaa
&&\frac{d^2 v}{d t^2}(t)=A v(t), \quad \frac{d^2 w}{d t^2}(t)=A w(t), \quad t\in (0,T),\label{eqR}\\
&&v(0)=\Re a,\quad w(0)=\Im a,\label{Ra0}\\
&&\int_0^T [\cos(\o t)v(t)-\sin(\o t)w(t)dt]=\Re g,\quad \int_0^T [\sin(\o t)v(t)+\cos (\o t)w(t)dt]=\Im g.\label{R}
\eaaa
 }
\end{remark}

\subsection*{Connection with the  Cauchy problem}
For $a\in\H$ and $b\in H$, consider a boundary  problem with the  Cauchy condition
\baa
&&\frac{d^2 u}{d t^2}(t)=A u(t), \quad t\in (0,T),\label{eq00}\\
&&u(0)=a, \label{a}
\\
&&\frac{d u}{d t}(0)=b\label{b}.
\eaa

For $u\in \W^{1}_1$, we accept again that equation (\ref{eq00}) is
satisfied  as an equality (\ref{intA}) that holds in $H^{-1}$ for all $s,t$ such that $0\le s\le t\le T$, and that conditions (\ref{a}) and (\ref{b}) are
satisfied as equalities  in $H^{-1}$.


\begin{proposition}
\label{prop1} 
For any  $a\in H^1$ and $b\in H$, there exists a unique solution
$u\in\W^1_1$ of problem (\ref{eq0})-(\ref{b}). This solution $u$ is such that
$u\in\W^{1}_\infty$. Moreover, there exists $c>0$  such
that
\baa \|u\|_{{\W}^1_\infty}\le c(
\|a\|_{\H}+\|b\|_{H})\quad
\forall a\in\H,\ b\in H.
\label{estp0} \eaa
Here $c>0$ depends only on $\HBC,A,T$, and $\o$.
\end{proposition}
The statement of Proposition \ref{prop1} represents a minor modification 
of well known results adapted to our choice of spaces; however, we provided its proof in Section \ref{SecP} below for the sake of completeness.

\section{Proofs}
\label{SecP}

{\em Proof of Proposition \ref{prop1}}.
 Let $a$  and $b$ be expanded  as
\baa
a=\sum_{k=1}^\infty \a_k v_k,\quad b =\sum_{k=1}^\infty \b_k v_k.
\label{ab}\eaa
Here the coefficients $\a_k$ are such that $\sum_{k=1}^{+\infty}|\a_k|^2=\|a\|_H<+\infty$.

We look for the solution $u$   expanded  as
\baa
u(t)=\sum_{k=1}^\infty  y_k(t)v_k,
\label{uy}
\eaa
where $y_k(t)$ are solutions of equations \baa
\frac{d^2 y_k}{d t^2} (t)=\bl{\rho \frac{d y_k}{d t} (t)}-\lambda_ky_k(t).
\label{y}
\eaa
In this case, \baaa \frac{d^2 u}{d t^2} (t)=\frac{d^2 }{d t^2}\sum_{k=1}^{\infty}y_k(t)v_k
=\sum_{k=1}^{\infty}(\bl{\rho \frac{d y_k}{d t} (t)}-\lambda_ky_k(t)v_k
=\sum_{k=1}^{\infty}y_k(t)Av_k=Au(t).\eaaa

Let \baaa \bl{p\defi \rho/2,\quad} \t_k\defi \sqrt{\bl{p^2+}\lambda_k}.
\eaaa
It can be seen that \baa y_k(t)=
C_ke^{\bl{pt}-i\t_kt}+D_ke^{\bl{pt+}i\t_k t}
\label{yCD}\eaa
for some $C_k,D_k\in\C$.
\bl{It can be noted that this representation holds for the case where $p^2+\lambda_k>0$
 as well as  for the case where $p^2+\lambda_k<0$; in this case, $\t_k$ is an imaginary number.
******If $p^2+\lambda_k=0$ - consider separately}

The coefficients  $C_k$ and $D_k$ are defined from the system
 \baa &&C_k+D_k=\a_k,\nonumber\\ && 
 -i\t_kC_k +i\t_kD_k=\b_k.
\bl{(p-i\t_k)C_k +(p+i\t_k)D_k=\b_k.}
\label{CD0}
\eaa
This gives
\baaa &&D_k=\a_k-C_k,\\ &&
(\a_k-D_k)(-i\t_k)  +D_ki\t_k=2i \t_kD_k  -\a_k i\t_k=\b_k. \eaaa
and \baaa
-\a_k i\t_k +i\t_k D_k
 +D_k i\t_k=2i\t_kD_k  +\a_k( - i\t_k)=\b_k.
\eaaa
Hence
\baa
D_k=\frac{\b_k+\a_ i\t_k}{2i\t_k},\quad  C_k=\frac{-\b_k+\a_k i\t_k}{2i\t_k}.
\label{CD00}
\eaa
                        
\bl{ \baaa
\a_k p-D_k p -\a_k i\t_k +i\t_k D_k
 +D_k p+D_k i\t_k=2i\t_kD_k  +\a_k(p - i\t_k)=\b_k.
\eaaa 
Hence
\baa
D_k=\frac{\b_k+\a_k(i\t_k-p)}{2i\t_k},\quad  C_k=\frac{-\b_k+\a_k(i\t_k+p)}{2i\t_k}.
\label{CD00}\eaa}
For the case of real $\a_k$ and $\b_k$, we have that $C_k=\oo D_k$.

For $p=0,1$ and $q=0,1$ such that $p+q=1$, we have that
 \baa
&&\sup_{t\in[0,T]}\left\|\frac{d^p u}{d t^p}(t)\right\|_{H^q}^2\le c_1\sum_{k=1}^{\infty}\lambda_k^{p+q}
|y_k(t)|^2
\le c_2\sum_{k=1}^{\infty}\lambda_k^{p+q} (|D_k|^2+|C_k|^2)\nonumber\\&&\le c_3
\sum_{k=1}^{\infty}\lambda_k^{p+q} \left(|\a_k|^2+\frac{|\b_k|^2}{{\lambda_k}}\right)
\le c_4(\|a\|_{H^{p+q}}+\|b\|_{H^{p+q-1}})=c_4(\|a\|_{H^{1}}+\|b\|_{H^0}).
\label{energy02}\eaa
This means that estimate (\ref{estp0}) holds.

Clearly, equations (\ref{eq0})--(\ref{b}) hold for the case where $a$, $b$, and $u$, are replaced by their  truncated expansions
\baa
a_N=\sum_{k=1}^N \a_k v_k,\quad b_N =\sum_{k=1}^N \b_k v_k,\quad u_N=\sum_{k=1}^N y_k(t) v_k.
\label{xiN}\eaa
We have that  $u_N\in \W^2\cap\W^1$. Estimate (\ref{energy02}) and completeness
of the Banach space $\W^1$ ensures that $u_N\to u$ in $\W^1$ as $N\to+\infty$. This $u$ is the solution (\ref{eq0})--(\ref{b}), and  that  energy estimate (\ref{energy02}) holds.
This completes the proof of Proposition \ref{prop1}. $\Box$

To proceed to the proof of Theorem \ref{ThM},
we need to adjust the approach used to  the case of the integral boundary conditions.

Let $a$  and $g$ be expanded  as
\baa
a=\sum_{k=1}^\infty \a_k v_k,\quad g =\sum_{k=1}^\infty \g_k v_k.
\label{ag_ser}\eaa
Here the coefficients $\a_k$ are such that $\sum_{k=1}^{+\infty}\lambda_k|\a_k|^2=\|a\|_\H<+\infty$.

We look for the solution $u$   expanded  as
\baa
u(t)=\sum_{k=1}^\infty  y_k(t)v_k,
\label{u_exp}
\eaa
where $y_k(t)$ are the solutions of equations (\ref{y})
defined by (\ref{yCD}), where
 $C_k$ and $D_k$ are defined from the system
 \baa &&C_k+D_k=\a_k,\nonumber\\ &&
 C_k \int_0^T e^{i\o t-it\t_k}dt  +D_k \int_0^T e^{i\o t+it\t_k}dt=\g_k.
\label{sys}\eaa
\begin{lemma}\label{lemmaCD}  Solution $(C_k,D_k)$ of system (\ref{sys}) exists and is uniquely defined for any $k$. Moreover,  there exists $c>0$ that depends on $T$ and $\o$ only
and such that
\baa
|C_k|+|D_k|\le c\left(|\a_k|+(1+\t_k)|\g_k|\right)
\label{CD}\eaa
for all $k$.
\end{lemma}

{\em Proof of Lemma \ref{lemmaCD}}.
Let 
\baaa
d_k \defi\int_0^T e^{i\o t-it\t_k}dt -\int_0^T e^{i\o t+it\t_k}dt.
\label{dkdef}
\eaaa

In this case,
\baaa&&C_k=\a_k-D_k,\nonumber\\ &&
\a_k\int_0^T e^{i\o t-it\t_k}dt  +D_k d_k=\g_k.
\label{Dg}
\eaaa

Suppose that we can prove that
\baa
 \d\defi \inf_{k} |d_k(1+\t_k)|>0.
\label{d}
\eaa
In this case, $d_k\neq 0$ for all $k$, and  there exists $c_p>0$ such that, for $k\in\Lambda_1^+$,  
\baaa
|D_k|&\le& c_1|d_k|^{-1} (2|\a_k/(1+\t_k)| +|\g_k|)=c_1|
d_k(1+\t_k)|^{-1} (2|\a_k| +(1+\t_k)|\g_k|)
\nonumber\\
&\le& {\d}^{-1}   (2|\a_k +(1+\t_k)|\g_k|).
\label{CD2}\eaaa
This would imply the proof of the lemma.

Let us prove that (\ref{d}) holds.
Let
\baa
&&\Lambda_0\defi \{k:\   
\hbox{either}\quad \t_k =\o  \quad \hbox{or}\quad \t_k =-\o\},\nonumber
\\
&&\Lambda_1\defi \{k\notin \Lambda_0:\ \hbox{either}\quad e^{i\t_k T}=e^{i\o T}\quad 
\hbox{or}\quad e^{i\t_k T}=e^{-i\o T}\},\nonumber\\
&&\Lambda_2\defi \{k:\ \quad e^{i\t_k T}\neq e^{i\o T},\quad e^{i\t_k T}\neq e^{-i\o T}\}.\label{Lambda}
\eaa
Clearly, these sets are disjoint, and the set $\Lambda_0$ is either an empty set or a singleton.

The assumption of Theorem \ref{ThM} that $e^{2i\o T}\neq 1$ excluded the case where
  $e^{i\t_k T}= e^{i\o T}$ and  $e^{i\t_k T}= e^{-i\o T}$  simultaneously, since this would imply that
 $e^{i\o T}= e^{-i\o T}$, or $e^{2i\o T}= 1$.
 Hence
 \baaa
 \Lambda_0\cup \Lambda_1\cup \Lambda_2=\{k=1,2,3,...\}.
 \eaaa
\subsubsection*{The case where $k\in\Lambda_0$}
Let us consider first the case where  $\Lambda_0=\{k\}$ is a singleton, 
  i.e., $\t_k=\o$ or $\t_k=-\o$.
  
Let us consider the case where  $\o=\t_k$.
In this case,
\baaa &&C_k=\a_k-D_k,\qquad (\a_k-D_k)T +D_k\w d=\g_k,
\eaaa
where
\baaa \w d=\frac{1}{i\o +i\t_k}\left[e^{i\o T+iT\t_k}-1\right]=
\frac{1}{2i\o}\left[e^{2i\o T}-1\right].
\label{dk1}\eaaa
Hence
$D_k (\w d-T)  =-\a_k T+\g_k.$ 
We have that  $\w d-T\neq 0$, since, clearly, 
$|e^{2it\o}|<|1+2iT\o|$. The case where  $\o=-\t_k$ can be considered similarly. 

\comm{ Hence
\baaa
\d_0\defi \inf_{k\in\Lambda_0} |d_k|>0.
\label{d}
\eaaa
Hence
\baaa
&&\sup_{k\in \Lambda_0} |D_k|\le \d_0^{-1} |\a_k T+\g_k|\le {\d_0}^{-1}(T|\a_k| +|\g_k|),\quad |C_k|\le |D_k|+|\a_k|.
\label{CD1}\eaaa
This implies that (\ref{CD}) holds for this singe $k$.
Therefore,
the statement of  Lemma \ref{lemmaCD} holds  if $k\in \Lambda_0$.}

\subsubsection*{The case where $k\in\Lambda_1$}

Let us consider the case where  $k\in \Lambda_1$. By the definitions, it follows that 
$\o\neq \pm\t_k$.
In this case, 
\baa d_k=
\frac{1}{i\o +i\t_k}\left[e^{i\o T+iT\t_k}-1\right]-\frac{1}{i\o -i\t_k}\left[e^{i T\o -i\t_k T}-1\right]
\label{dk2}
\eaa
and
\baa
D_k d_k  =-\frac{\a_k\left[e^{i\o T-iT\t_k}-1\right]}{i\o -i\t_k} +\g_k.
\label{Dd}
\eaa

\par
Let us consider first  the case where  $e^{i\t_k T}=e^{i\o T}$.
In this case, we have that
\baaa d_k=
\frac{1}{i\o+i\t_k }\left[e^{2i\o T}-1\right].
\label{dk3}
\eaaa
By the assumptions on $\o$ and $T$, $e^{2i\o T}\neq 1$.
 Hence
 \baaa
 \inf_{k\in\Lambda_1^+} |d_k|(1+\t_k)>0,
\label{d1+}
\eaaa 
where $\Lambda_1^+=\{k:\ e^{i\t_k T}= e^{-i\o T}\}$.

 Using a similar approach where $e^{i\t_k T}= e^{-i\o T}$, we obtain that 
\baaa
\inf_{k\in\Lambda_1} |d_k|(1+\t_k)>0,
\label{d1}
\eaaa
Hence the statement of Lemma \ref{lemmaCD} holds for $k\in\Lambda_1$.

\subsubsection*{The case where $k\in\Lambda_2$}
Let us consider  the most typical  case where  $k\in \Lambda_2$, i.e.
\baa
e^{i\t_k T}\neq e^{i\o T}, \quad e^{i\t_k T}\neq e^{-i\o T}.
\label{ee}\eaa
In particular,  we have in this case that $\o\neq \t_k$ and $\o\neq -\t_k$ for all $k$ in
this case.

Let us show first  that the values $d_k(i\o+i\t_k)$ are separated from zero for large $k$.
We have that
\baaa  d_k(i\o +i\t_k )&=&e^{i\o T -iT\t_k}\biggl[e^{2iT\t_k}-\frac{\o+\t_k}{\o -\t_k}\biggr]
-\frac{2\t_k}{\t_k-\o}\\&=&e^{i\o T}\biggl[e^{iT\t_k}+e^{ -iT\t_k}\biggr]+\xi_k
-\frac{2\t_k}{\t_k-\o}=2e^{i\o T}\cos(2T\t_k)+\xi_k
-\frac{2\t_k}{\t_k-\o},
\eaaa
where  \baaa
\xi_k\defi  e^{i\o T -iT\t_k}\biggl[-1-\frac{\o+\t_k}{\o -\t_k}\biggr]
=e^{i\o T -iT\t_k}\biggl[\frac{\t_k+\o}{\t_k -\o}-1\biggr]
=e^{i\o T -iT\t_k}\frac{2\o}{\t_k -\o}.
\eaaa
\comm{comm: Hence
\baaa  d_k(i\o +i\t_k)=2e^{i\o}\cos(2T\t_k)+e^{i\o T-iT\t_k}\frac{2\o}{\t_k -\o}
-\frac{2\t_k}{\t_k-\o}
\eaaa}
Hence
\baaa  \frac{d_k}{2}(i\o+i\t_k)=e^{i\o T}\Biggl[\cos(2T\t_k)+\zeta_k\Biggr]-z_k,\eaaa
where \baaa
z_k=\frac{\t_k}{\t_k-\o},
 \quad \zeta_k=\frac{1}{2}e^{ -i\o T}\xi_k.
\eaaa
Clearly, we have that    $\zeta_k\to 0$ and $z_k\to 1$ as $k\to +\infty$. Let  $a_k=\Re \zeta_k$ and $b_k=\Im \zeta_k$.
We have that $\zeta_k=a_k+ib_k$ and \baaa
\lim_{k\to +\infty}(a_k\cos(\o T) -b_k\sin (\o T) -z_k)=-1.\eaaa
Further, we have that
\baaa  \Re\frac{d_k}{2}(i\o+i\t_k)
=\cos (\o T) \cos(2T\t_k)+a_k\cos(\o T) -b_k\sin (\o T) -z_k.
\eaaa
Since $e^{2i\o T}\neq 1$, it follows that  $|\cos (\o T)|<1$ and  $\sup_k \cos (\o T)\cos(2T\t_k)<1$.
It follows that there exists $N>0$ such that
\baaa\sup_{k\ge N}\Re\frac{d_k}{2}(i\o+i\t_k)<0.
\eaaa
Hence
\baa
\inf_{k\ge N}(|\o|+|\t_k|)|d_k|>0.
\label{N}
\eaa

To complete the proof, it suffices to show that  $d_k\neq 0$  for $\t_k\in\Lambda_2$.
We have that
\baaa  d_k=\frac{2f(\t_k)}{i(\o^2-\t_k^2)},
\eaaa
where
\baaa  f(x)=\frac{1}{2}[(\o-x)(e^{i\o T+ix T}-1) -(\o+x)(e^{i\o T-ixT}-1)]\\
=\frac{1}{2}\Bigl(e^{i\o T} [\o (e^{i x T}-e^{-i x T}) -  x(e^{ix T}+e^{-ix T}) ]+2x\Bigr)\\
=e^{i\o T} [i\o \sin ( x T)-  x\cos(x T) ]+x.
\eaaa
Let us show that $f(x)\neq 0$ for all $x\in\{\t_k\}_{k\in\Lambda_2}$.
Suppose that $f(x)=0$. In this case,
\baa
&&\Im f(x)=\o\cos(\o T)\sin (xT)-x \sin(\o T)\cos(x T)=0,
\nonumber
\\ && \Re f(x)=-\o\sin(\o T)\sin (xT)-x \cos(\o T)\cos(x T)+x=0.
\label{IR}\eaa

Suppose that equations (\ref{IR}) hold and that $\cos(xT)=0$. In this case,
$\o\cos(\o T)\sin (xT)=0$ and $\sin (xT)\neq 0$. This implies that $\cos(\o T)=0$
and that either $e^{i\o T}=e^{i xT}$ or $e^{i\o T}=-e^{i xT}$.
Similarly, suppose that equations (\ref{IR}) hold and that $\cos(\o T)=0$. In this case,
$x\sin(\o T)\cos (xT)=0$ and $\sin (\o T)\neq 0$. This implies  that $\cos(\o T)=0$.
Again, it follows  either $e^{i\o T}=e^{i xT}$ or $e^{i\o T}=-e^{i xT}$.
Hence $d_k\neq 0$ in both cases, since, as
is shown above,   $d_k\neq 0$ if  $e^{i\o T}=e^{i \t_k T}$ or $e^{i\o T}=-e^{i \t_k T}$
for some $k$.

Therefore, it suffices to consider the case where $k\in\Lambda_2$  and $\cos(x T)\neq 0$ and $\cos(\o T)\neq 0$.

The equation for $\Im f$  in (\ref{IR}) gives that
\baaa
\o\frac{\sin (xT)}{\cos (xT)}=x\frac{\sin(\o T)}{\cos (\o T)},
\eaaa
i.e.
\baaa
\tan(xT)=\frac{x}{\o}\tan (\o T),\quad  \sin (xT)= a(\o)x\cos(xT),
\eaaa
where $a(\o)\defi\o^{-1}\tan (\o T)$.
Then the equation for $\Re f$ in (\ref{IR})   gives that
\baaa
-\o\sin(\o T) a(\o)x\cos(xT)-x \cos(\o T)\cos(x T)+x=0.
\eaaa
This can be rewritten as
\baaa
\sin(\o T)\tan (\o T)x\cos(xT) +x \cos(\o T)\cos(x T)-x=0.
\eaaa
Hence
\baaa
\cos(x T)x(\sin^2(\o T) +\cos^2(\o T))/\cos (\o T)-x=0
\eaaa
and
\baaa
\cos(x T)x /\cos (\o T)-x=0.
\eaaa
This implies that
\baaa
x(1-\cos(x T)/\cos (\o T))=0.
\eaaa
This would imply that $\cos(x T)=\cos (\o T)$ and $\Re e^{i x T}=\Re e^{i \o T}$, and this in turn would imply  that either $e^{i x T}=e^{i \o T}$
or $e^{i x T}=e^{-i \o T}$.  However, this case
is excluded for $k\in \Lambda_2$.
Therefore, $f(\t_k)\neq 0$ and $d_k\neq 0$   for $k\in \Lambda_2$. Hence
\baa
\inf_{k\in\Lambda_2} |d_k(1+\t_k)|>0.
\label{N1}
\eaa
Hence  (\ref{d}) holds.
This  completes the proof of Lemma \ref{lemmaCD}. $\Box$

We are now in the position to prove Theorem \ref{ThM}.

{\em Proof of Theorem  \ref{ThM}}. Let $u$
be defined by (\ref{u_exp}),(\ref{yCD}),(\ref{sys}).  We have that, for the case of $\g$ and $a$
 this $u$
 is a unique solution  problem of  (\ref{eq})--(\ref{ppp}) as well as
 problem  (\ref{eq0})--(\ref{b}) with  $b=du(0)/dt$; this can be shown
similarly to the
proof of Proposition \ref{prop1}. The identity is straightforward for truncated eigenfunction expansions with finite number of terms.  On the next step,  the extension on the case of infinite expansions is ensured by the energy estimates.

 We have that
 \baaa
 \frac{du}{dt}(0)=\sum_{k=1}^\infty (-i\t_kC_k+ i\t_kD_k )v_k=
 \sum_{k=1}^\infty  (-i\t_k(\a_k-D_k)+ 2i\t_k D_k)v_k
=\sum_{k=1}^\infty i\t_k (-\a_k+2   D_k )v_k.
 \eaaa
By Lemma \ref{lemmaCD}, it follows that $D_k^2\le c_0 (\a_k^2+\lambda_k \g_k^2)$
for some $c_0>0$ that depends on $\HBC$, $A$, $T$, and $\o$.  Hence
 \baa
\left\| \frac{du}{dt}(0)\right\|_{H}^2
\le c_1 \sum_{k=1}^\infty \lambda_k  (|C_k|^2+ |D_k|^2)\le c_2 \sum_{k=1}^\infty \lambda_k
[\a_k^2+
 (1+\t_k)^2 \g_k^2)]\le
  c_3(\|a\|_{\H}^2+\|g\|_{H^2}^2)
\eaa
for some $c_1>0$, $c_2>0$,  and $c_2>0$,  that depend on $\HBC$, $A$, $T$, and $\o$.
By Lemma \ref{prop1},  estimate (\ref{estp}) holds.
 This
 completes the proof of Theorem \ref{ThM}. $\Box$

\section{A numerical example of impact of the presence of $\o\ne 0$}
\label{SecN}
Consider a toy example for the problem
 \baaa && u''_{tt}(x,t)=u''_{xx}(x,t),\qquad t\in[0,T],\
x\in(0,\pi)\\ &&u(x,0)=0,\quad 
\int_0^T e^{i\o}u(x,t)dx=g(x), \quad
u(0,t)=u(\pi,t)=0,
\eaaa
where $g\in H^2$. This is a special case of problem (\ref{eq})-(\ref{ppp}) with $n=1$, $D=(0,\pi)$, $A=d^2/d x^2$.   
It is known that $\lambda_k=k^2$ and 
$v_k(x)=\sin(kx)$,  $k=1,2,...,$,  are the corresponding eigenvalues and eigenfunctions
Respectively, $\t_k=k$.  \bl{see PDEtrent, example 7.7 label{exevef1}}

For simplicity, we assume that $\o\notin\{1,2,3,...\}$.  In this case, $\t_k\in\Lambda_2$,
 in the notations of the  proof of Theorem \ref{ThM}. In addition, we assume that $e^{2i\o T}\neq 1$. 
The solution constructed in the proof of Theorem \ref{ThM} is  defined by (\ref{dk2})
and (\ref{De}), i.e.,
\baaa
u(x,t)=\sum_{k=0}^{+\infty}y_k(t)u_k(x), 
\eaaa 
where \baaa y_k(t)=
C_ke^{-i\sqrt{\lambda_k}t}+D_ke^{i\sqrt{\lambda_k}t}, 
\eaaa
and where $C_k$ and $D_k$ are defined from the system (\ref{dk2}), (\ref{Dd}) such that
\baa
C_k=-D_k, \quad D_k  =\frac{\g_k}{d_k}.
\label{De}
\eaa
where $d_k$ defined by (\ref{Dd}),  i.e.
\baaa d_k=
\frac{1}{i\o +i k}\left[e^{i\o T+iTk}-1\right]-\frac{1}{i\o -ik}\left[e^{i T\o -i k T}-1\right].
\label{dk20}
\eaaa
Let $z(m)\defi \inf_{k\le m}|d_k|(1+\t_k)$. As can be seen from the proof of Lemma \ref{lemmaCD}, this values  should be separated from zero to ensure regularity of solutions claimed in Theorem \ref{ThM}. 

If $\o=0$, then the problem of small denominators arises: for any $T>0$, $z(m)\to 0$ as $m\to +\infty$ and hence $|D_k|\to \infty$. 

 In particular, for $T=5$, we estimated  numerically  that  $z(500)=3.66\cdot 10^{-9}$ for  $\o=0$ and that $z(500)=0.1001$ 
for $\o=0.01$.   For $T=10$, we estimated  numerically  that 
$z(500)=3.68\cdot 10^{-9}$ for  $\o=0$ and that $z(500)=0.1998$ 
for $\o=0.01$. 

This illustrates that that the impact of including even a small enough $\o\neq $
can prevent appearance of "small denominators" (small divisors).

\section{Conclusions and  discussion}\label{SecD}
The paper establishes solvability and regularity of a complexified boundary value problem
for linear hyperbolic wave  equations  where a Cauchy condition is replaced by
 a  integral condition $\int_0^Te^{i\o t}u(t)dt=g$ for the solution.
It is shown that this new problem is well-posed in a wide class of solutions
given the presence of  the weight function $e^{i\o t}$, where $\o\in\R\setminus\{0\}$ can be arbitrarily small; the only condition is that $e^{2iT\o}\neq 1$.   This leads  to complex valued solutions of the boundary value problem with this integral condition. This boundary value problem would be equivalent to a boundary value problem 
for a system of two real valued
 hyperbolic equations, for the real and imaginary parts of the complex valued solution
 respectively, with
an integral condition connecting solutions. In this case, the real part of the solution can be considered as an approximation as $\o\to 0$ of the solution of the real valued solution with $\o=0$. 
The setting considered in the paper allows many modifications and extensions.
Most likely, the results can be extended on the case
where the eigenvalues for $A$ can be non-positive, and where the weight
$e^{i\o t}$ in (\ref{ppp}) is replaced by $e^{rt+i\o t}$ for $r\in\R$. We leave it for the future research.

So far, we have considered the case where the solutions can be expanded via the basis from the eigenfunctions.  It would be interestingly to extend the result  on
the more general case, as it was done in \cite{GG} for wave equations with two point conditions.
We leave it for the future research as well.
\subsubsection*{Acknowledgment}
This work was supported in part by the Zhejiang University/University of Illinois at Urbana-Champaign Institute.

 \end{document}